# The new axiom of set theory and Bell inequality.

## Kajetan Guz

Extension and intension are connected with the distinguishability of objects. It is important whether we are dealing with distinguishable or indistinguishable entities, and how they coexist among themselves. Is distinguishability a feature that belongs to mathematical entities independently of cognitive acts, or are those features perhaps assigned by consciousness in the course of creative analysis? This problem was considered by Zermelo when constructing the axiom of separation *(Axiom der Aussonderung)*. In his work of 1906 there appears the term "definite", used to express the fact that all elements of a domain are pairwise distinguishable in such a way that in every case it is possible to ascertain whether a=b or a≠b. Zermelo refers to "well-defined sets" *(wohldefinierte Mengen)* and "well-defined properties" *(wohldefinierte Eigenschaften)*. The term "definite" was used and analysed by Edmund Husserl in two lectures given at the Göttingen Mathematical Society. A deeper mathematical analysis of that concept appeared in his "Formal and Transcendental Logic" in relation to a definite aggregate [*1*]. We shall return to that work later. Bertrand Russell, analysing Plato's arguments contained in the "Parmenides" dialogue [*Introduction to Mathematical Philosophy*] comes to the conclusion that numbers do not have being, and calls them "logical fictions". Stewart Shapiro considered the Julius Caesar problem in a chapter titled "Structuralism–Ontology: Object" [*2*]. It is a problem that orients the philosophy of mathematics towards Ontology.
The problem arises again in analysis by Benacerraf addressing the issue of multiple reduction, namely the representation of numbers in different ways [3]. What kind of set is the number 3 – {∅, {∅},{∅,{∅}}} in von Neumann's scheme, or {{{∅}}} in Zermelo's – or is it a set at all? Benacerraf takes the view that we do not perceive numbers as sets.
Stewart Shapiro says that the philosophy of mathematics does not need to answer questions of the type "does Gajus Julius Cesar = 3?" or "is 1ϵ 3?" It ought rather to show why those questions do not require an answer.

   2. Primary concepts

Going back to early Greek thought, we find the word λεγω "collect", "select", "accumulate". This verb is the origin of the noun λόγος, which from Homer and Hesiod through Heraclitus and Plato took a permanent place in Greek philosophy. Both λεγω and λέγειν also mean "speaking". It may be wondered what connects speaking with collecting, this being a question about the connection of a distributive with a collective function. To quote Krzysztof Narecki: "The decisive factor here was an innovation in the form of a figurative meaning of λέγειν: 'proceed, (follow) a list, enumerate [cf. Latin *legere enatum*], relate, narrate'. Just as spatially (physically) one goes through (looks through) a certain number of objects in order to collect them in a single place and as a single whole, that same act may be performed discursively, that is, merely by the use of reason one may go through (track) a series of words to combine them into a single, logical (i.e. meaningful) statement. 'Speaking' (like 'collecting') thus becomes a synthetic action, arising from the etymological meaning of λέγειν" [4].

Georg Cantor, inspired by the ideas of Plato, initiated a new groundbreaking theory – the Set theory [5]. The name "set theory" *(Mengenlehre)* from very early on gained primacy over the original term "theory of multiplicities" *(Mannigfaltigkeitslehre)* and came to be used universally in all mathematical publications. Cantor defines an multiplicity *(Mannigfaltigkeit)* in parallel with the concept of a set *(Menge)* as any multiplicity *(jedes Viele)* which can be considered as a unity *(Eines)*, a generality of defined elements which by dint of some law may be unified into a single whole. He distinguishes a definite multiplicity *(Vielheit)*, which can be considered as a kind
of whole, unity, as a "ready thing" *(ein fertiges Ding)*, which he calls an uncontradictory multiplicity or set, from a definite multiplicity which does not fulfil that condition and which he calls an absolutely infinite or contradictory multiplicity. Until the end of his life, Cantor would struggle with the task of giving a philosophical justification of his theory based on mathematics.

The need for a philosophical foundation for the new theory was noted by Edmund Husserl, especially since he was both a mathematician and philosopher. In his "Logical Researches" [6] he explains the ideas of pure set theory. That work would enable Husserl, many years later, to make a profound philosophical analysis of the foundations of mathematics in *Formale und transzendentale Logik*. During those years he observed attentively the development of the new ideas, which inspired German mathematicians in their studies of the new mathematical entities. He was able to carry on a dispute with the mathematician-philosophers (criticism of Frege's *Die Grundlagen der Arithmetik* in *Philosophie der Arithmetik*, and Frege's reply in his review of that work), draw conclusions from it, and in some matters concede that his opponents were right. Significant and important analysis of the fundamental concepts was provided by Józef M. Bocheński [7]. Bocheński introduced two theories: the theory of similarity and the theory of identity.

The theory of similarity classifies things according to their similarity. The relation of similarity is a symmetric and reflexive relation. It concerns the bond existing between things being classified.

The theory of identity concerns the bond of certain properties that pertain in all objects being classified. It is a symmetric, reflexive and transitive relation. These theories imply that we have to deal with two types of entities with different ontological status. For things to exist and for properties to exist are entirely different matters. Quine's criterion – to exist is to be the value of a variable – is not applicable in this case. Willard Van Orman Quine attached great importance to this problem, although he formulated it in different categories [8]. He bases his analysis on the traditional distinction between general terms and abstract singular terms ("square" and "squareness"). He alludes to Frege's principle of substitution of equals by equals. Abstract terms may appear in identity-related contexts, being subject to the laws of identity. Quine notes that general along with singular terms are treated equally, as names of separate objects. Detailed analysis of intension and extension in the axioms of set theory brings to light what might henceforth be called the problem of the empty set. Many authors adopt various systems for the set of axioms of set theory, giving them the common name ZFC. Few, however, permit themselves to make a philosophical analysis of the content of each axiom. And such an analysis is necessary. Frege was not mistaken when he claimed: "A mathematician who is not a philosopher is only half a mathematician."

3. Axioms

Let us consider the axiom of extensionality, which is listed first in all of the systems of axioms.
If sets a and b have the same elements, then they are identical.

In symbols:

I. $\forall a \forall b (\forall x (x \in a \Leftrightarrow x \in b) \Rightarrow a = b)$

The axiom of extensionality comes from Ernst Zermelo, who formulated it as the *Axiom der Bestimmtheit* ("axiom of decisiveness") in 1907: *"Ist jedes Element einer Menge M gleichzeitig Element von N und umgekehrt, ist also gleichzeitig M ϵ N und N ϵ M, so ist immer M=N* [9]. The first formulation of the axiom was made by Bolzano [*Grössenlehre*], but was not published until after his death, in 1975 [10], but this had no influence on Zermelo's work. Independently of Bolzano, the axiom was formulated by Richard Dedekind in an essay of 1888 [11] and it was on that formulation that Zermelo's was based. Now, in all works on set theory, the term "axiom of extensionality" or "axiom of extension" is used.
The axiom of extensionality is not an existential axiom. It does not postulate the existence of a set. It is a criterion for the equality of two given sets whose elements can be compared. Commentaries on this axiom highlight the legitimacy of extension, and at the same time indicate the possibility of forming intensional concepts. The roots of this axiom go back to the 17th century, to the principle of identity of indiscernibles *(principium identitatis indiscernibilium)* of Gottfried Wilhelm Leibniz, which can be presented in the form:

II. $\forall x \forall y [\forall P (Px \Leftrightarrow Py) \Rightarrow x = y]$

For any x and y, if x and y have the same properties, then x is identical to y.

This law is schematically the same as the axiom of extensionality. It is not *principium magnum, grande et nobilissimum* as Leibniz described the *principium rationis* – a basic principle for all representation, where any being is regarded as existing if and only if it has been secured for representation as a calculable object (Martin Heidegger). The law of identity of indiscernibles is nonetheless controversial in its essence (Max Black): objects may be presented having the same properties but not being identical (numerically nonidentical). Leibniz also believed that there cannot exist two things that are perfectly essentially similar. He tried to derive that principle from the Principle of Reason *(nihil est sine ratione)* [12]. The Principle of Reason directs a person to calculate. Reason is understood as *ratio,* as a calculation. A compact philosophical analysis of the principle was given by Martin Heidegger [*13*]. It is possible that Zermelo called his first axiom the "axiom of decisiveness" because for Leibniz the principle of identity of indiscernibles was a principle of individualisation. Extension was meant to determine a set definitively as a set. Quine accurately stated that we need sets as extensional distillates of properties [*14*]. Verifying that the elements of sets are the same requires reference to their properties – that is, to intension. Such verification or measurement of properties is nonetheless a step in the direction of their distinguishability. The question arises of whether *wohldefinierte Eigenschaften* well-define a set. When the set is a quantum system, measurement of a property (distinguishability) changes the set. The Hong–Ou–Mandel experiment [15] involving the interference of two photons, as well as many other experiments in quantum physics, show that measurement of the properties of elements of a set has an effect on the entire set.

There is no controversy about the converse rule (the indiscernibility of identicals):

III. $\forall x \forall y [x = y \Rightarrow \forall P (Px \Leftrightarrow Py)]$

For any x and y, if x is identical to y, then x and y have the same properties.

Let us consider a set consisting of eight elements: six line segments of unit length 1, and two segments of length √2. This can be pictured in the following system:

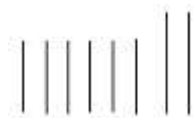

Fig. 1

or in another system:

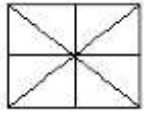
Fig. 2

The latter is called Quine's square, because he considered such a system when analysing the problem of universals [*8*].
In both examples there appear identical objects. By performing reduction, in the first system we obtain two objects from eight, while in the second, from 33 objects, we obtain five (a right-angled isosceles triangle, a square, a rectangle and two forms of general quadrilateral).
In his analysis, Quine went further, demonstrating the possibility of reducing the five shapes to one square within whose area they appear.
The question arises as to whether the two sets in Fig. 1 and Fig. 2 are equal. Starting from the axiom of extensionality we find that the sets are identical, because they have the same elements. If we analyse the sets taking an intensional approach, we obtain two different sets. The set in Fig. 1 is a distributive set, while that in Fig. 2 is a collective set. It might be added that the set in Fig. 1 is an actually distributive and potentially collective set, while that in Fig. 2 is actually collective and potentially distributive. The need to pay attention to and to distinguish distributive from collective properties was pointed out by Edmund Husserl in his *Philosophie der Arithmetik,* rejecting Frege's charge that he had not taken into account that in the phrase *vier edel Rose*, *vier* and *edel* are properties of different degree. *Vier* is not a distributive property stated of every rose, but a property stated of a collective whole. If we permit geometric thinking, we might transform the figure in Fig. 1 into that of Fig. 2, which in turn proves to be a set of different figures. To make it possible to compare the two sets according to the axiom of extensionality, the elements of those sets must be distinguishable. Such an assumption is important in set theory, since it enables proof of the existence of a "scale of alephs". If the elements of the two sets were indistinguishable, then we would have to count the elements in each of them to check their extension. Counting assigns a feature to each element, and the set of those features (numerosity) makes it possible to compare the sets. The axiom of extensionality is a recipe for the "measurement" of sets. The measure is the counting of the elements of the set.
In the axiom of extensionality, the relation "is an element of", denoted "$\epsilon$", is reflexive and symmetric, while the equality relation "=" is reflexive, symmetric and transitive. If we take the relation "$\epsilon$" to be transitive also, then we are dealing with a collective understanding of a set. The relation of equality has an intensional basis, while the "is an element of" relation has an extensional basis. This difference means that this axiom, being the criterion for the equality of two sets, is also a weighing scale for intensional and extensional concepts.

The second of the axioms on Zermelo's list – although omitted from many contemporary workings of set theory, ranging from textbooks [16] to would-be philosophical analyses [17] – is the empty set axiom.

There exists a set which does not have any element

Symbolically:

$$\text{IV.} \quad \exists a \forall b \quad \neg(a \in b)$$

From the start this axiom aroused much controversy. Zermelo calls this set improper *(uneigentliche Menge)*. Bertrand Russell posed the question of whether a non-empty set remains a set if we remove all of its elements [*18*]. In a letter to Abraham Fraenkel in 1921, Zermelo expressed doubts as to the justification for the empty set. He even suggested giving up that axiom, with the help of an appropriate restriction of the axiom schema of specification/separation.
There are also works which, with certain modifications to set theory (the system introduced by von Neumann and analysed by Robinson, Bernays and Gödel – NBG, VNB), accept this axiom [19]. A compact analysis of the empty set in the historical context is given by Akihiro Kanamori [*20*]. In essence this is a detailed discussion of axiom II on Zermelo's list, which was called the axiom of elementary sets *(Axiom der Elementarmengen)*. For this reason Kanamori considers the empty set, singleton and ordered pair within a single topic. In his remarks on the second axiom [*21*] Zermelo states, referring to the first axiom, that there exists only one empty set, while the singleton and ordered pair are unique. In a further remark [*22*] he states that the empty set is a subset of any set.
The empty set axiom is the most important axiom, alongside the axiom of infinity, among the axioms of set theory. These two are existential axioms – they refer to the existence of mathematical entities. All of the others are "constructive" axioms, which we shall call C axioms (the subset axioms, pair axiom, axiom of union, power set axiom, and axioms of replacement, regularity and choice). The C axioms allow one to deduce the existence of constructed sets from the assumption of the existence of certain other sets.
From Zermelo onwards it has become accepted to deduce from the axiom of extensionality that there is only one empty set (in many identical copies, as
a subset of every set). This reasoning is false, as it is based on a categorical error. The empty set axiom, together with the axiom of infinity, is superordinate to the other axioms. The empty set has an empty actual extension, but it also has an infinite potential intension.
From the axiom of extensionality one can deduce the uniqueness of the sets derived from the C axioms, except for the axiom of choice. One can prove the uniqueness of a power set; but in attempting to prove the uniqueness of the empty set, that scheme is inappropriate.

### 3. Empty sets

The problem of the empty set in mathematics, and failure to recognise its significance, has its roots in the beginnings of European philosophical-mathematical thought, and reached a zenith in Leibniz's Principle of Reason. Emptiness did not have any reason, any foundation, any attribute to exist in the ideal world of mathematical entities. The first philosopher to consider this problem was Edmund Husserl [1]. This work, completed in 1929 after more than a decade of silence, was written "(…) in one go and printed in just a few months – after I had thought for decades over the issues taken up" (from a letter to Roman Ingarden). It was a period when the development of set theory was taken up by mathematicians who were not interested in seeking philosophical foundations for it. Only Husserl's student Roman Ingarden, also a philosopher and mathematician, appreciated the great potential of the aforementioned work. In *Formale und transzendentale Logik* Husserl defines the concepts of "set" and "cardinal number" at a level of pure and broadest generality, referring to "objectiveness in general", "something in general" – in the most empty generality. This formality requires any material definition to be something indefinitely arbitrary. The basic concepts of set theory and the theory of cardinal numbers are syntactic creations, syntactic forms of the empty **Something**. Husserl sees in the derivative forms of something in general the unity of formal mathematics, which in all of its related disciplines has one basis – the empty region as something in general.

Let us consider one more axiom, the axiom of union *(Axiom der Vereinigung)*, the fifth on Zermelo's list.

Symbolically:

V. $\forall A \exists B \forall c (c \in B \iff \exists D (c \in D \land D \in A))$.

For any set *A* there exists a set **U** *A* which contains exactly the elements of the elements of that set, and is called the union of *A*. From the first axiom we know that there can be only one such set. In this case that deduction does not produce doubts as in the case of the empty set; nonetheless, the validity of the axiom of extensionality, as we have previously analysed, does give rise to such doubts. These considerations are important in view of the application of unions of sets in quantum physics.

Let us introduce into our further deliberations a new axiom:

The axiom of existence of sets:

There exist infinitely many essentially different empty sets. Symbolically:

VI. $\forall \emptyset_a \forall \emptyset_b (\emptyset_a \neq \emptyset_b)$

An important property of empty sets is their uniqueness. If two sets are equal,

then their empty sets are different. Symbolically:

$$\text{VII.} \quad \forall a \forall b (a = b \Rightarrow \emptyset_a \neq \emptyset_b)$$

In contrast to equal sets, sets that are identical to each other have identical empty sets.
Symbolically:

$$\text{VIII.} \quad \forall a (a = a \Rightarrow \emptyset_a = \emptyset_a)$$

What is at issue here is the Principle of Identity. The essence of every set is its empty set. Important deliberations on the Principle of Identity were undertaken by Martin Heidegger [23]. Heidegger sees this Principle as a form of a fundamental principle which takes identity as a characteristic in being, that is, in the ground of a being *(als einen Zug im Sein, d. h. im Grund des Seienden voraussetzt)*. Taking a leap away from this Principle, abandoning being as the ground of a being, one moves to groundlessness *(Abgrund)*. However, this groundlessness is neither an empty nothingness or murky chaos, but an occurrence, denoted as *Er-eignis (Doch dieser Abgrund ist weder das leere Nichts noch eine finstere Wirrnis, sondern: das Er-eignis)*. The empty set is the groundlessness for a set, it is not an empty nothingness but an occurrence for the set: it shapes and creates the set.

An important property of empty sets is their probabilistic nature, which is revealed in interactions with other empty sets. Every set contains within itself an empty set with a probability dependent on the probabilities of the empty sets belonging to that set.

Symbolically:

$$\text{IX.} \quad \forall a (\exists \emptyset_a \Rightarrow P\emptyset_a = \cup \emptyset_A)$$

This probability is the sum of the probabilities of the empty sets of the elements of that set. We now apply axiom V to empty sets:

$$\text{X.} \forall A \exists B \forall \emptyset_c (\emptyset_c \in B \Leftrightarrow \exists D(\emptyset_c \in D \land D \in A))$$

In this case set *B* does not form the union of the empty sets of set *A*. Let us consider a simple example:

$$A = \{\{a,b\},\{b,c\},\{a,c\}\}$$
$$\cup A = \{a,b,c\}$$

but: $\cup \emptyset_A = P\emptyset_{ab} + P\emptyset_{bc} + P\emptyset_{ac}$

Another property of empty sets is the differentiation of the cardinality of sets by empty sets.
Symbolically:

$$\forall A \forall B(|A| \neq |B| \Leftrightarrow P\emptyset_A \neq P\emptyset_B)$$

If any two sets have different cardinalities, then the probabilities of their empty sets are different – and conversely, if the probabilities of any two empty sets are different, then the cardinalities of the sets are different.

## 4. Quantum reality

Quantum theory, which has enjoyed much success in physics, has raised questions to which it is difficult to give precise answers. One of these is the interpretation of Bell's Theorem (Bell's inequality). John Bell stated his theorem in 1964 [24]. It was presented in a simple way, in the language of set theory, using Venn diagrams, by Leonard Susskind (Fig. 3).

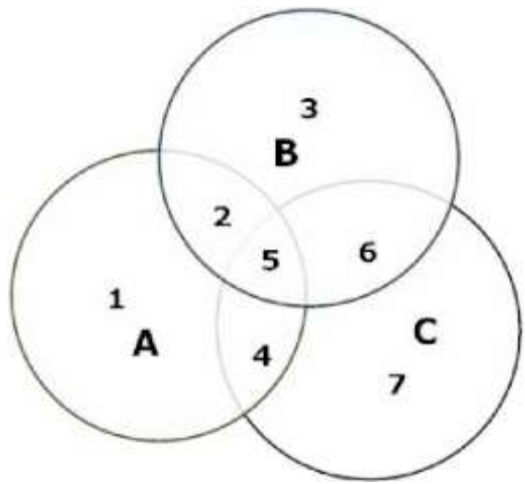

Fig. 3

According to the logic of set theory, we can write the following inequality:

>A not B + B not C ≥ A not C

or

>(1+4) + (3+2) ≥ (1+2)
>(1+2) + (4+3) ≥ (1+2)

which holds.
Bell's inequality is written in the form of an inequality of probabilities:

>P(A not B) + P(B not C) ≥ P(A not C)

However, experiments in quantum physics contradict this inequality. All interpretations to date are directed against the thought experiment of Einstein,

Podolsky and Rosen (EPR).

Attention has not been paid, however, to the imperfection of the mathematical apparatus used to describe quantum reality. Using the new axiom of empty sets we obtain a different Venn diagram, which can be used to present Bell's inequality in a different form.

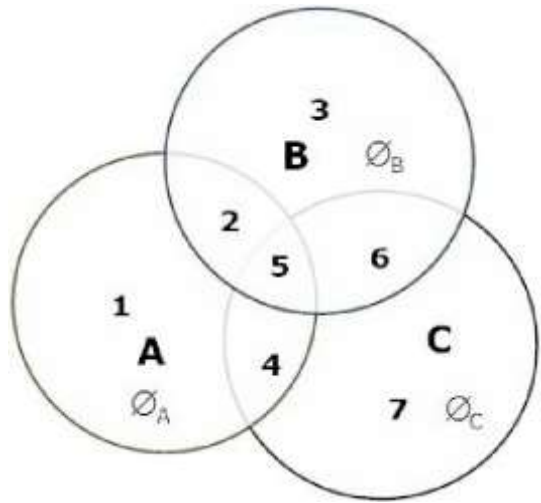

Fig. 4

Each of the three sets A, B and C has its respective empty set: $\emptyset_A$, $\emptyset_B$ and $\emptyset_C$. The following inequality holds:

$$P\emptyset_{A\bar{B}} + P\emptyset_{B\bar{C}} \neq P\emptyset_{A\bar{C}}$$

Bell's inequality takes the form:

P(A not B) + P(B not C) ≠ P(A not C)

Bibliography

1. Edmund Husserl , *Formale und transzendentale Logik. Versuch einer Kritik der logischen Vernunft.* Max Niemeyer Verlag, Halle (Saale)